\documentclass[12pt]{article}

\def\stackunder#1#2{\mathrel{\mathop{#2}\limits_{#1}}}

\begin{document}
\title{q-Analogue of $A_{m-1}\oplus A_{n-1}\subset A_{mn-1}$}
\author{V. G. Gueorguiev, A. I. Georgieva, P. P. Raychev and R. P. Roussev\\
\\
\it Institute for Nuclear Research and Nuclear Energy,\\
\it Bulgarian Academy of Science\\ 
\it 1784 Sofia,Bulgaria}

\date{PACS. 02. 20 - Group theory. }
\maketitle

\begin{abstract}
A natural embedding $A_{m-1}\oplus A_{n-1}\subset A_{mn-1}$ for the
corresponding quantum algebras is constructed through the appropriate
comultiplication on the generators of each of the $A_{m-1}$ and $A_{n-1}$
algebras. The above embedding is proved in their $q$-boson realization by
means of the isomorphism between the $\mathcal{A}_q^{-}$ (mn)$\sim {\otimes}
^n \mathcal{A}_q^{-}$(m)$\sim {\otimes}^m\mathcal{A}_q^{-}$(n) algebras. 
\end{abstract}

Recently, quite a great interest has been paid to the study of quantum
algebras and their applications to physical problems. Essentially quantum
algebras are Hopf algebras. Hopf algebra is an algebra with additional
structures: i) except the multiplication $m:A\otimes A\rightarrow A$ there
is a comultiplication $\Delta :A\rightarrow A\otimes A$; ii) except the unit 
$1$ which provides the embedding $R\rightarrow A$ ($C\rightarrow A$), where $
R$ ($C$) is the real (complex) field, there is a counit $\varepsilon
:A\rightarrow R\;(C)$ . All these mappings are homomorphisms and there is an
antihomomorphism $S:A\rightarrow A$ called antipode. From mathematical point
of view such algebras were developed much earlier \cite{HA1,HA2}. The
contemporary development of their theory is connected with noncommutative
geometry and differential calculus \cite{geometry}. In physics these new
mathematical objects appear in the theory of the inverse scattering problem 
\cite{ISP}. Later on, the quantum algebras have been applied to a number of
physical regions as statistical mechanics , quantum field theory , molecular
, atomic and nuclear physics. 

In nuclear structure theory successful applications of models, based on
algebraic chains of Lie algebras ( Interaction Boson Model (IBM) \cite{IBM}
, Two Vector Boson Model (TVBM) \cite{TVBM} etc. ) were obtained. 

It will be of interest to construct q-analogues of these chains and consider
the corollary of the models based on them. The chain $su_q(3)\oplus
u_q(2)\subset u_q(6)$ was already considered by Quesne in \cite{Quesne}. 

In this paper we consider the general case of the embedding : 
\begin{equation} \label{emb}
A_{m-1}^q\oplus A_{n-1}^q\subset A_{mn-1}^q
\end{equation}
in q-boson realization. The proper definition of the embedding (\ref{emb})
is a result of a careful analysis of the comultiplication structure. The
present paper also provides the method of its realization, briefly described
bellow. 

As, it is well known, for any integer $n$ the algebra $A_{n-1}^q$ has a
realization of its generators in terms of the q-boson algebra $\mathcal{A}
_q^{-}(n)$ \cite{q-boson1,q-boson2}. In order to obtain the realization of
the generators of $A_{n-1}^q$ in terms of the $q$-boson algebra $\mathcal{A}
_q^{-}$(mn), we apply $m-1$ times the comultiplication, then a q-boson
realization for each term in the tensor product and finally employ the
isomorphism $\mathcal{A}_q^{-}$(mn) $\sim {\otimes}^m\mathcal{A}_q^{-}$(n). 
By analogy we realize the generators of $A_{m-1}^q$. The generators of the $q
$ -deformed algebra $A_{mn-1}^q$ have their realization by means of the same
algebra $\mathcal{A}_q^{-}$(mn) . 

We start with the algebraic relations among the regular functionals $
l_{ij}^{\pm}$ of the quantum matrix group given in \cite{Reshetikhin}:

\begin{equation} 
\label{eq:RLL}
\begin{array}{cc}
{\stackunder{m,p}{\sum}R_{ij,mp}^{+}l_{mk}^{+}l_{pl}^{-}}={\stackunder{m,p}{
\sum}l_{jp}^{-}l_{im}^{+}R_{mp,kl}^{+}}\;; & {\stackunder{m,p}{\sum}
R_{ij,mp}^{+}l_{mk}^{\pm}l_{pl}^{\pm}}=
{\stackunder{m,p}{\sum}l_{jp}^{\pm}l_{im}^{\pm}R_{mp,kl}^{+}}
\end{array}
\end{equation}

In the case of deformed $A_{n-1}^q$ algebras the explicit form of the $R^{+}$
-matrix is given by: 
\begin{equation} \label{eq:R}
R^{+}=q^{\frac 1n}\{q\sum_{i=1}^ne_{ii}\otimes e_{ii}+\sum_{i\ne
j=1}^ne_{ii}\otimes e_{jj}+(q-q^{-1})\sum_{i<j=1}^ne_{ij}\otimes e_{ji}\}
\end{equation}
where $e_{ij}$ are $n\times n$ matrixes with elements $(e_{ij})_{km}=\delta
_{ik}\delta _{jm}$. 

By substituting (\ref{eq:R}) in (\ref{eq:RLL}) we obtain the following
relations for $l_{ij}^{\pm}$ : 
\begin{equation}
\begin{array}{c}
\lbrack l_{im}^{(\varepsilon )},l_{js}^{(\varepsilon )}]=(1-q)( \underbrace{
l_{im}^{(\varepsilon )}l_{js}^{(\varepsilon )}}_{i=j}- \underbrace{
l_{js}^{(\varepsilon )}l_{im}^{(\varepsilon )}} _{m=s})+(q-q^{-1})(
\underbrace{l_{jm}^{(\varepsilon )}l_{is}^{(\varepsilon )}}_{m>s}-
\underbrace{l_{jm}^{(\varepsilon )}l_{is}^{(\varepsilon )}}_{j>i}) \\ 
{\lbrack l_{im}^{+},l_{js}^{-}]}=(1-q)(\underbrace{l_{im}^{+}l_{js}^{-}}
_{i=j}-\underbrace{l_{js}^{-}l_{im}^{+}}_{m=s})+(q-q^{-1})(\underbrace{
l_{jm}^{-}l_{is}^{+}}_{m>s}-\underbrace{l_{jm}^{+}l_{is}^{-}}_{j>i}) \\ 
{\ \stackrel{n}{\stackunder{i=1}{\prod}}}l_{ii}^{\pm}=1\quad ;\quad 
l_{ii}^{+}l_{ii}^{-}=1=l_{ii}^{-}l_{ii}^{+}\quad \\
\quad {l_{ij}^{+}=0} 
\quad \mathrm{for}\quad i>j\quad \mathrm{and}\quad l_{ij}^{-}=0\quad 
\mathrm{for}\quad i<j
\end{array}
\end{equation}

The last relations employs not only the form of $R^{+}$ but also some
additional conditions \cite{Reshetikhin}. 

Further by means of the substitution: 
\begin{equation} \label{eq. Yl}
\begin{array}{c}
H_{ij}=\tilde H_i-\tilde H_j\; \\ 
\;l_{ij}^{\pm}=\mp q^{\pm {\ \frac 12}}(q-q^{-1})Y_{ij}^{\pm}q^{\mp {\frac
12}(\tilde H_i+\tilde H_j)}\;\mathrm{with}\;Y_{ii}^{\pm}={\mp}{\frac{
q^{\mp {\ \frac 12}}}{{q-q^{-1}}}}
\end{array}
\end{equation}

one comes to the following relations for the Cartan-Weyl basis of q-deformed 
$A_{n-1}^q$ algebra: 
\[
\begin{array}{c}
\lbrack H_{ij},H_{km}]=0 \\ 
\begin{array}{ll}
\quad Borel \quad subalgebra:\mathcal{B}^{+} & \quad Borel \quad
subalgebra:\mathcal{B}^{-} \\ 
{\lbrack Y_{ik}^{+},Y_{kj}^{+}]}_q=Y_{ij}^{+}\quad i<k<j & {\
[Y_{ij}^{-},Y_{jk}^{-}]}_{q^{-1}}=Y_{ik}^{-}\quad i>j>k \\ 
{\lbrack Y_{ik}^{+},Y_{ij}^{+}]}_q=0\quad i<j<k & {[Y_{kj}^{-},Y_{ij}^{-}]}
_{q^{-1}}=0\quad i>k>j \\ 
{\lbrack Y_{kj}^{+},Y_{ij}^{+}]}_q=0\quad i<k<j & 
{[Y_{ik}^{-},Y_{ij}^{-}]}_{q^{-1}}=0\quad i>j>k \\ 
{\lbrack Y_{ij}^{+},Y_{km}^{+}]}=0\quad i<j<k<m & 
[Y_{ij}^{-},Y_{km}^{-}]=0 \quad i>j>k>m \\ 
{\lbrack Y_{ij}^{+},Y_{km}^{+}]}=0i<k<m<j & 
[Y_{ij}^{-},Y_{km}^{-}]=0 \quad i>k>m>j \\ 
{\lbrack Y_{km}^{+},Y_{ij}^{+}]} =(q-q^{-1})Y_{kj}^{+}Y_{im}^{+} & 
[Y_{ij}^{-},Y_{km}^{-}]=(q-q^{-1})Y_{kj}^{-}Y_{im}^{-} \\ 
i<k<j<m & i>k>j>m \\ 
{\lbrack H_{ik},Y_{js}^{+}]}=(e_i-e_k,e_j-e_s)Y_{js}^{+} & {\
[H_{ik},Y_{js}^{-}]}=(e_i-e_k,e_j-e_s)Y_{js}^{-}
\end{array}
\end{array}
\]

\begin{equation} \label{eq:algebraic}
\begin{array}{c}
Mixed \quad commutators \\ 
\\ 
{\lbrack Y_{ij}^{+},Y_{ji}^{-}]}=[H_{ij}]_qi<j
\end{array}
\end{equation}
\[
\begin{array}{ll}
{\lbrack Y_{km}^{-},Y_{ij}^{+}]}=(q-q^{-1})Y_{kj}^{+}Y_{im}^{-}q^{H_{ik}} & 
{\ [Y_{ij}^{+},Y_{km}^{-}]}=(q-q^{-1})Y_{kj}^{-}Y_{im}^{+}q^{H_{jm}} \\ 
j>k>i>m & k>j>m>i \\ 
{\lbrack Y_{ij}^{+},Y_{im}^{-}]}=0j>i>m & {[Y_{ij}^{+},Y_{kj}^{-}]}=0k>j>i
\\ 
{\lbrack Y_{ij}^{+},Y_{ki}^{-}]}=-Y_{kj}^{+}q^{H_{ik}}j>k>i & {\
[Y_{ij}^{+},Y_{ki}^{-}]}=-q^{H_{ji}}Y_{kj}^{-}k>j>i \\ 
{\lbrack Y_{ij}^{+},Y_{jm}^{-}]}=Y_{im}^{-}q^{H_{ij}}j>i>m & {
[Y_{ij}^{+},Y_{jm}^{-}]} =q^{H_{jm}}Y_{im}^{+}j>m>i \\ 
{\lbrack Y_{ij}^{+},Y_{km}^{-}]}=0 & \left\{ 
\begin{array}{l}
k>j>i>m;\;k>m>j>i \\ 
j>k>m>i;\;j>i>k>m
\end{array}
\right. 
\end{array}
\]
where $(e_i,e_j)={\delta}_{ij}$ , the $q$-commutator is given by ${[A,B]}
_q=AB-qBA$ and the $q$-number is defined by $[x]_q={\frac{{q^x-q^{-x}}}{{\
q-q^{-1}}}}$. These relations are analogical to the ones obtained in \cite
{Birow}. 

It should be noted here that the generators $Y_{ij}^{\pm}$ can be
substituted by $\tilde Y_{ij}^{\pm}f_{ij}(q,\tilde H)$, which will lead to
modifications in the relations (\ref{eq:algebraic}) depending on the
functions $f_{ij}(q,\tilde H)$. An example of such a mapping from $su(2)$ to
a deformed $su_q(2)$ is given in \cite{Zachos}. 

From the definition of the \textit{Comultiplication} $\Delta (l_{ij}^{\pm})=
\sum_{k=1}^nl_{ik}^{\pm}\otimes l_{kj}^{\pm}\;$ and the \textit{Counit} 
$\;\varepsilon (l_{ij}^{\pm})=\delta _{ij}$ given in \cite{Reshetikhin} we
obtain the following coalgebraic structure : 
\begin{equation} \label{eq:coalgebraic}
\begin{array}{l}
\Delta H_{ij}=H_{ij}\otimes 1+1\otimes H_{ij}\;;\;\varepsilon
(H_{ij})=0\;;\;S(H_{ij})=-H_{ij} \\ 
\varepsilon (Y_{ij}^{\pm})= {\mp}{\frac{q^{{\mp}{\frac 12}}}{{q-q^{-1}}}}{
\delta}_{ij}\;;\;{\ Y_{ii}^{\pm}}={\mp}{\frac{q^{{\mp}{\frac 12}}}{{
q-q^{-1}}}}\;;\;{\ Y_{ik}^{+}}=0\;i>k\;;\;{Y_{ik}^{-}}=0\;i<k \\ 
{\Delta Y_{ij}^{\pm}}={\mp} (q-q^{-1})q^{{\pm}{\frac 12}}\displaystyle{
\sum_{i\le k\le j\;\mathrm{or} \;(j\le k\le i)}}Y_{ik}^{\pm}q^{{\pm}{\frac
12}{H_{jk}}}\otimes Y_{kj}^{\pm}q^{{\pm}{\frac 12}{H_{ik}}}
\end{array}
\end{equation}

Applying the standard definition of the antipode $S$ ($m\circ ({id}\otimes
S)\circ{\Delta} = m\circ (S\otimes id)\circ {\Delta} = i\circ {\varepsilon}$
) we deduce for the antipode of the generators $Y_{ij}^{\pm}$ the following
recurrent formula: 
\begin{equation}
S(Y_{ij}^{\pm})= -q^{\mp 1}Y_{ij}^{\pm}\;{\pm}(q-q^{-1})q^{\pm 1} 
\displaystyle{\sum_{i<k<j \; \mathrm{or}\; (i>k>j)}} Y_{ik}^{\pm}S(
Y_{kj}^{\pm})
\end{equation}

Let us introduce the q-boson algebra $\mathcal{A}_q^{-}$(n) with creation
and annihilation operators $a_i^{\pm}$ and their q-boson numbers $N_i$ as
in \cite{q-boson1,q-boson2,q-boson3,q-boson4}. 
\begin{equation} \label{eq:boson}
a_i^{-}a_i^{+}-q^{\mp}a_i^{+}a_i^{-}=q^{\pm N_i}\;\mathrm{and}
\;[N_i,a_j^{\pm}]=\pm {\delta}_{ij}a_j^{\pm}
\end{equation}

The $q$-boson realization of the Cartan-Chevalley generators $H_i=H_{i,i+1}$
, $Y_i^{+}=Y_{i,i+1}^{+}$ and $Y_i^{-}=Y_{i+1,i}^{-}$ of $A_{n-1}^q$-algebra
given by Sun and Fu in \cite{q-boson1} is: 
\begin{equation} \label{eq:Chevalley}
H_i=N_i-N_{i+1}\;;\;Y_i^{+}=a_i^{+}a_{i+1}^{-}\;;
\;Y_i^{-}=a_{i+1}^{+}a_i^{-}
\end{equation}

The irreducible Fock representations $\Gamma _q^{[m]}$ with the vacuum state 
$\mid 0>$, ${b_i^{-}|0>=0}$, $N_i\mid 0>=0$ is defined by the set of
vectors: 
\begin{equation} \label{Fock irrep}
\Gamma _q^{[m]}:=\{\mid m>=\mid m_1,. . . ,m_n>=\displaystyle{\ \prod_{i=1}^n{
\frac{(b_i^{+})^{m_i}}{{\sqrt{[m_i]!}}}}}\mid 0>\;\;\mid m={\ \sum_{i=1}^nm_i}\}
\end{equation}
with the following properties: 
\begin{equation} \label{Fock2}
\begin{array}{l}
dim\Gamma _q^{[m]}= {\frac{(n+m-1)!}{m!(n-1)!}} \\ 
N\mid m>=m\mid m>\;\;\mathrm{where}\;\;N= \displaystyle{\sum_{i=1}^nN_i}. 
\end{array}
\end{equation}

Using the definitions of $H_i$ (\ref{eq:Chevalley}) and $N$ (\ref{Fock2})
the operators $N_i$ can be expressed by 
\begin{equation} \label{N_i by N and H}
N_i={\frac 1n}N+{\frac 1n}\displaystyle{\
\sum_{s=2}^n\sum_{j=1}^{s-1}H_j-\sum_{j=1}^{i-1}H_j}
\end{equation}

The additional generators which extend (\ref{eq:Chevalley}) to the basis ( 
\ref{eq:algebraic}) of Cartan-Weyl can be obtained from the Chevalley
generators (\ref{eq:Chevalley}) by means of the first relations in the Borel
subalgebras $\mathcal{B}^{\pm}$ in (\ref{eq:algebraic}). In this way, as in 
\cite{Quesne1} we obtain the following general realization: 
\begin{equation} \label{eq:Ya}
\begin{array}{lr}
H_{ij}=N_i-N_j\;;\; & Y_{ij}^{\pm}=a_i^{+}a_j^{-}q^{{\mp}\displaystyle{\
\sum_{i<k<j\;\mathrm{or}\;(j<k<i)}}N_k}
\end{array}
\end{equation}

Let us denote the generators of $A_{k_1k_2-1}^q$ by $Y_i^{\pm}$ and $N_i$,
of $A_{k_1-1}^q$ by $X_{\mu}^{\pm}$ and $N_{\mu}$, of $A_{k_2-1}^q$ by $
Z^{\pm s}$ and $N^s$ and the $n$-th product of the comultiplication by: 
\[
{\Delta}^n=(\underbrace{id\otimes id\otimes \ldots \otimes \Delta}_n)( 
\underbrace{id\otimes id\otimes \ldots \otimes \Delta}_{n-1})\ldots
(id\otimes \Delta)\Delta 
\]

Since $\Delta$ is a homomorphism one can consider the following mapping: 
\begin{equation} \label{eq:homomorphism}
A_{m-1}^q\; \displaystyle{\rightarrow^{\Delta^{(n-1)}}} \; \underbrace
{A_{m-1}^q \otimes \ldots \otimes A_{m-1}^q}_{n}
\end{equation}

For the sake of simplicity, the tensor product $\otimes $ will be dropped
and the index $s$ $($or $\mu )$ will indicate the number of the tensor
space. Thus we obtain: 
\begin{equation} \label{eq:delta}
\begin{array}{ll}
{\tilde H}_\mu ={\stackrel{k_2}{\stackunder{s=1}{\sum}}}H_\mu ^s; & {\tilde
X}_\mu ^{\pm}=\Delta ^{(k_2-1)}(X_\mu ^{\pm})={\stackrel{k_2}{\stackunder{
s=1}{\sum}}}X_\mu ^{\pm s}q^{{\frac 12}\displaystyle{\sum_{\sigma \ne
s,\sigma =1}^{k_2}}sign(\sigma -s)H_\mu ^\sigma} \\ 
{\tilde H}^s={\stackrel{ k_1}{\stackunder{\mu =1}{\sum}}}H_\mu ^s; & {
\tilde Z}^{\pm s}=\Delta ^{(k_1-1)}(Z^{\pm s})={\stackrel{k_1}{\stackunder{
\mu =1}{\sum}}}Z_\mu ^{\pm s}q^{{\frac 12}\displaystyle{\sum_{\sigma \ne
\mu ,\sigma =1}^{k_1}} sign(\sigma -\mu )H_\sigma ^s}
\end{array}
\end{equation}
\bigskip

From the construction of the operators (\ref{eq:delta}) and as a result of
the used homomorphism $\Delta $ it is easy to prove that the generators $
\tilde X_\mu ^{\pm}$, $\tilde H_\mu $ and $\tilde Z^{\pm s}$, $\tilde H^s$
satisfy the commutations relations for the algebras $A_{k_1-1}^q$ and $
A_{k_2-1}^q$. 

Using the q-boson realization of the generators (\ref{eq:Ya}) we obtain : 
\begin{equation} \label{eq:Ak1}
\begin{array}{l}
{\tilde X}_\mu ^{+}=\displaystyle{\sum_{s=1}^{k_2}}a_\mu ^{+s}a_{\mu
+1}^{-s}q^{{\frac 12}\displaystyle{\sum_{\sigma \ne s,\sigma =1}^{k_2}}
sign(\sigma -s)(N_\mu ^\sigma -N_{\mu +1}^\sigma )} \\ 
{\tilde X}_\mu ^{-}= \displaystyle{\sum_{s=1}^{k_2}}a_{\mu +1}^{+s}a_\mu
^{-s}q^{{\frac 12} \displaystyle{\sum_{\sigma \ne s,\sigma =1}^{k_2}}
sign(\sigma -s)(N_\mu ^\sigma -N_{\mu +1}^\sigma )} \\ 
{\tilde Z}^{+s}=\displaystyle {\sum_{\mu =1}^{k_1}}a_\mu ^{+s}a_\mu
^{-s+1}q^{{\frac 12}\displaystyle{\sum_{\sigma \ne \mu ,\sigma =1}^{k_1}}
sign(\sigma -\mu )(N_\sigma ^s-N_\sigma ^{s+1})} \\ 
{\tilde Z}^{-s}=\displaystyle{\sum_{\mu =1}^{k_1}}a_\mu ^{+s+1}a_\mu ^{-s}q^{
{\frac 12}\displaystyle{\sum_{\sigma \ne \mu ,\sigma =1}^{k_1}} sign(\sigma
-\mu )(N_\sigma ^s-N_\sigma ^{s+1})} \\ 
{\tilde H}^s= \displaystyle{\sum_{\mu =1}^{k_1}}{N_\mu ^s-N_\mu ^{s+1}}\;;\;{
\tilde H}_\mu =\displaystyle{\sum_{s=1}^{k_2}}{N_\mu ^s-N_{\mu +1}^s}
\end{array}
\end{equation}

It is correct to consider the q-bosons in $\tilde X$ and $\tilde Z$ (\ref
{eq:Ak1}) as different objects, because in $\tilde X$, $a_\mu ^{\pm s}$ mean:

$a_\mu ^{\pm s}=\underbrace{id\otimes \ldots \otimes id\otimes \overbrace{
a_\mu ^{\pm}}^s\otimes id\otimes \ldots \otimes id}_{k_2}$

\noindent while in $\tilde Z$ :

$a_\mu ^{\pm s}=\underbrace{id\otimes \ldots \otimes id\otimes \overbrace{
a_s^{\pm}}^\mu \otimes id\otimes \ldots \otimes id}_{k_1}$

\noindent However in both cases, they satisfy the same relations: 
\begin{equation} \label{eq:relboson}
\begin{array}{ll}
\lbrack a_\mu ^{\pm s},a_\nu ^{\pm t}]=0\;\mathrm{for\ all}\;s,t,\mu ,\nu & 
[a_\mu ^{+s},a_\nu ^{-t}]=0\; \mathrm{for\ all}\;s\ne t;\mu \ne \nu \\ 
{\lbrack N_\mu ^s,a_\nu ^{\pm t}]} =\pm \delta _{\mu ,\nu}\delta
_{s,t}a_\nu ^{\pm t} & a_\mu ^{-s}a_\mu ^{+s}-q^{\mp 1}a_\mu ^{+s}a_\mu
^{-s}=q^{\pm N_\mu ^s}
\end{array}
\end{equation}
Let us define the following correspondence: $i\leftrightarrow (\mu ,s)$ ($
k_2\le k_1$): 
\begin{equation} \label{isomorphic}
\begin{array}{l}
i\leftrightarrow (\mu ,s)\;\;i=1,\ldots ,k_1k_2;\;\mu =1,\ldots
,k_1;\;s=1,\ldots ,k_2 \\ 
\mu =1+ \mathrm{int}{[\frac{i-1}{k_2}]}\;\;\mathrm{where}\;\mathrm{int}
\left[ x\right] \;\mathrm{is}\;\mathrm{integer}\;\mathrm{part}\;\mathrm{of}
\;x \\ 
s=1+(i-1)\mathrm{mod} (k_2)\;,\;i=(\mu -1)k_2+s
\end{array}
\end{equation}
\bigskip

From the introduction of (\ref{isomorphic}) in equations (\ref{eq:boson})
and (\ref{eq:relboson}) it follows that The algebras ${\otimes}^{k_2}
\mathcal{A}_q^{-}(k_1)$ and ${\otimes}^{k_1}\mathcal{A}_q^{-}(k_2)$
constructed by the $q $-bosons $a_\mu ^{\pm s}$ are isomorphic to the
algebra $\mathcal{A} _q^{-}(k_1k_2)$ constructed by the $q$-bosons $a_i^{\pm}$. 
As a result the algebras $A_{k_1-1}^q$ and $A_{k_2-1}^q$ have
realization in the $\mathcal{A} _q^{-}$($k_1k_2$) algebra. \bigskip

\textbf{Proposition 1. } \textit{The the generators ${\tilde X}_\mu ^{\pm}$
, ${\ \tilde H}_\mu $ commute with the generators ${\tilde Z}^{\pm s}$ , ${
\tilde H}^s$ given by (\ref{eq:Ak1}) . }\bigskip

{\it Proof. }
Let us consider the commutator between the elements ${\tilde X}_\mu ^{+}$
and ${\tilde Z}^{-s}$. For this purpose we define $Q_{t,\nu}$ and 
$I_{t,\nu}(\mu ,s,k)$ as: 
\[
Q_{t,\nu}=q^{{\frac 12}\displaystyle{(\sum_{\sigma \ne t,\sigma
=1}^{k_2}sign(\sigma -t)(N_\mu ^\sigma -N_{\mu +1}^\sigma )+\sum_{\rho \ne
\nu ,\rho =1}^{k_1}sign(\rho -\nu )(N_\rho ^s-N_\rho ^{s+1}))}} 
\]
\[
I_{t,\nu}(\mu ,s,k,q)=q^{{\frac 12}\displaystyle{\sum_{\sigma \ne t,\sigma
=1}^ksign(\sigma -t)(\delta _{\mu ,\nu}-\delta _{\mu +1,\nu})(\delta
_{\sigma ,s+1}-\delta _{\sigma ,s})}} 
\]
Using (\ref{eq:Ak1}) and (\ref{eq:relboson}), for the commutator we obtain: 
\begin{equation} \label{eq:proof}
{[{\tilde X}_\mu ^{+},{\tilde Z}^{-s}]=\sum_{t=1,\nu =1}^{k_2,k_1}}\{a_\mu
^{+t}a_{\mu +1}^ta_\nu ^{+s+1}a_\nu ^sI_{t,\nu}(\mu ,s,k_2,q)-a_\nu
^{+s+1}a_\nu ^sa_\mu ^{+t}a_{\mu +1}^tI_{\nu ,t}(s,\mu
,k_1,q^{-1})\}Q_{t,\nu}
\end{equation}
The sum over $t$ and $\nu $ can be represented as a sum of five terms: 
\[
\begin{array}{ll}
(\mathrm{a})=\{\nu \ne \mu ,\mu +1\;\mathrm{and}\;t\ne s,s+1\} & ( \mathrm{b}
)=\{\nu =\mu \;\mathrm{and}\;t=s+1\} \\ 
(\mathrm{c})=\{\nu =\mu +1\;\mathrm{and}\;t=s\} & ( \mathrm{d})=\{\nu =\mu
\;\mathrm{and}\;t=s\} \\ 
(\mathrm{e})=\{\nu =\mu +1\;\mathrm{and} \;t=s+1\} & 
\end{array}
\]
In these cases we have: 
\[
I_{t,\nu}(\mu ,s,k_2,q)=\left\{ 
\begin{array}{ll}
1 & \mathrm{in\ (a)} \\ 
q^{\frac 12} & \mathrm{in\ (b),(d)} \\ 
q^{-{\frac 12}} & \mathrm{in\ (c),(e)}
\end{array}
\right. \;I_{\nu ,t}(s,\mu ,k_1,q^{-1})=\left\{ 
\begin{array}{ll}
1 & \mathrm{in\ (a)} \\ 
q^{\frac 12} & \mathrm{in\ (b),(e)} \\ 
q^{-{\frac 12}} & \mathrm{in\ (c),(d)}
\end{array}
\right. 
\]

In the cases (a), (b) and (c) the bosons $a_\nu ^{+s+1}$, $a_\nu ^s$, $a_\mu
^{+t}$ and $a_{\mu +1}^t$ commute and the relevant terms are equal to zero. 
Thus the commutator is given only by the sum of (d) and (e) i. e. 
\[
{[{\tilde X}_\mu ^{+},{\tilde Z}^{-s}]=}q^{-{\frac 12}}a_\mu ^{+s+1}a_{\mu
+1}^s(q^{-N_{\mu +1}^{s+1}}Q_{s+1,\mu +1}-q^{-N_\mu ^s}Q_{s,\mu})=0 
\]

The expression $sign(\rho -\mu )=sign(\rho -\mu -1)$ when $\rho <\mu $ or $
\rho >\mu +1$ is used essentially in the calculation of $q^{-N_{\mu
+1}^{s+1}}Q_{s+1,\mu +1}=q^{-N_\mu ^s}Q_{s,\mu}$ . The other commutators
can be proved in the same way

Further using (\ref{eq:Ya}) and the isomorphism (\ref{isomorphic}) we have: 
\begin{equation} \label{aaY}
a_{1+\mathrm{int}\left[ {\frac{i-1}{k_2}}\right]}^{1+(i-1)\mathrm{mod}
(k_2)}a_{1+\mathrm{int}\left[ {\frac{j-1}{k_2}}\right]}^{1+(j-1)\mathrm{mod}
(k_2)}=a_i^{+}a_j^{-}=Y_{ij}^{\pm}q^{\pm \displaystyle{\sum_{i<\sigma <j\; 
\mathrm{or}\;(i>\sigma >j)}}N_\sigma}
\end{equation}

Finally applying (\ref{N_i by N and H}) and (\ref{aaY}) the generators of $
A_{k_1-1}^q$ and $A_{k_2-1}^q$ in (\ref{eq:Ak1}) are expressed through the
generators of $A_{k_1k_2-1}^q$ in the following way: 
\begin{equation} \label{eq:Z}
\begin{array}{l}
{\tilde Z}^{\pm s}=\displaystyle{\sum_{\mu =1}^{k_1}}
Y_{(\mu -1)k_2+s}^{\pm}q^{{\frac 12}
\displaystyle{\ \sum_{\sigma \ne \mu ,\sigma
=1}^{k_1}} sign(\sigma -\mu )H_{(\sigma -1)k_2+s}} \\ 
{\tilde H}^s=\displaystyle{\ \sum_{\mu =1}^{k_1}}H_{(\mu -1)k_2+s}\;\;;\;\;{
\tilde H}_\mu =\displaystyle{\ \ \sum_{s=(\mu -1)k_2+1}^{(\mu -1)k_2+k_2}}
H_{s,s+k_2} \\ 
{\tilde X}_\mu ^{+}=\displaystyle{\sum_{t=\mu k_2+1}^{(\mu +1)k_2}}
Y_{t-k_2,t}^{+}q^{{\frac 12\ \displaystyle\sum_{\nu \ne t,\nu =\mu
k_2+1}^{(\mu +1)k_2}sign(\nu -t)H_{\nu -k_2,\nu}+\Lambda _t^{+}}} \\ 
{\tilde X}_\mu ^{-}=\displaystyle{\ \sum_{t=\mu k_2+1}^{(\mu +1)k_2}}
Y_{t,t-k_2}^{-}q^{{\ \frac 12\ \displaystyle \sum_{\nu \ne t,\nu =\mu
k_2+1}^{(\mu +1)k_2}sign(\nu -t)H_{\nu -k_2,\nu}+\Lambda _t^{-}}} \\ 
\Lambda _t^{\pm}={\frac{k_2-1}{{k_1k_2}}(N+ \displaystyle\sum_{\sigma
=2}^{k_1k_2}H_{1,\sigma})\pm \displaystyle \sum_{\sigma
=t-k_2+1}^{t-1}H_{1,\sigma}}
\end{array}
\end{equation}

The difference $\Lambda _t^{\pm}$ between the expressions for \textit{\ $
\tilde Z^{\pm s}$} and \textit{$\tilde X_\mu ^{\pm}$} is due to the
ordering of indices in (\ref{isomorphic}) which leads to the appearance of
different terms $q^{{\mp}\displaystyle{\sum_{i<k<j\;\mathrm{or}\;(j<k<i)}}
N_k}$ in the q-boson realization (\ref{eq:Ya}) of the Chevalley and the
additional Weyl generators. In the expression $\Lambda _t^{\pm}$ the
operator $N$ , in q-boson realization has the meaning of a total number of
bosons operator. In general a corresponding operator may be constructed in
some extension of the algebra \textit{$A_{k_1k_2-1}^q$}. This can be proved
by induction. For \textit{$A_1^q$} ($su_q(2)$) the operator $N$ can be
obtained from the second order Casimir operator:

\[
C_2^q=X^{-}X^{+}+[H/2]_q[H/2+1]_q={\frac{q^{N+1}+q^{-N-1}-q-q^{-1}}{
(q-q^{-1})^2}} 
\]
For $n>2$ , $N^{(n)}$ -- the corresponding operator $N$ for \textit{$
A_{n-1}^q$} , is obtained by the recurrence:

\begin{equation} \label{N(s) by N(s-1) and H}
{N^{(n)}}={\frac{n+1}n}\{N^{(n-1)}+{\frac
1{n+1}\sum_{t=2}^{n+1}\sum_{p=1}^{t-1}H_p-\sum_{p=1}^nH_p\}}
\end{equation}
Moreover in practice it is only the eigenvalues of $q^N$ which are required. 
\bigskip\ 

\textbf{Proposition 2. } \textit{The elements $\tilde X_\mu ^{\pm}$, $\tilde
H_\mu $ of $A_{k_1-1}^q$ and $\tilde Z^{\pm s}$, $\tilde H^s$ of $A_{k_2-1}^q
$ defined by (\ref{eq:Z}) belong to the algebra $A_{k_1k_2-1}^q$ and provide
an explicit embedding $A_{k_1-1}^q\oplus A_{k_2-1}^q\subset A_{k_1k_2-1}^q$
in the $q$-boson realization (\ref{eq:Ya}) of $A_{k_1k_2-1}^q$. } \bigskip

{\it Proof. }
From the above it follows that the elements defined by (\ref{eq:Z}) belong
to the $q$-deformed $A_{k_1k_2-1}^q$ algebra. Applying the q-boson
realization (\ref{eq:Ya}), the correspondence (\ref{isomorphic}) and (\ref
{eq:relboson}) we obtain the $q$-boson realization ( \ref{eq:Ak1}) of the
generators \textit{$\tilde X_\mu ^{\pm}$, $\tilde H_\mu $} and \textit{$
\tilde Z^{\pm s}$, $\tilde H^s$}, whose commutation relations close the
algebras \textit{\ $A_{k_1-1}^q$} and \textit{$A_{k_2-1}^q$}. Finally these
two pairs of generators commute between themselves as proved in \textbf{
Proposition 1}and so they close the algebra \textit{$A_{k_1-1}^q\oplus
A_{k_2-1}^q$} embedded in \textit{$A_{k_1k_2-1}^q\diamondsuit $}

The results of Quesne \cite{Quesne} are reproduced in the case $k_1k_2=6$ , $
k_1=3$ and $k_2=2$. 

In the limit $q\rightarrow 1$ we obtain the usual embedding: 
\[
\begin{array}{lll}
{\tilde H}_\mu ={\stackrel{(\mu -1)k_2+k_2}{\stackunder{s=(\mu -1)k_2+1}{
\sum}}}H_{s,s+k_2} & {\tilde X}_\mu ^{+}={\stackrel{k_2}{\stackunder{s=1}{
\sum}}}Y_{(\mu -1)k_2+s,\mu k_2+s}^{+} & {\tilde X}_\mu ^{-}=
{\stackrel{k_2}{\stackunder{s=1}{\sum}}}Y_{\mu k_2+s,(\mu -1)k_2+s}^{-} \\ 
{\tilde H}^s={\ \stackrel{k_1}{\stackunder{\mu =1}{\sum}}}H_{(\mu -1)k_2+s}
& {\tilde Z} ^{\pm s}={\stackrel{k_1}{\stackunder{\mu =1}{\sum}}}Y_{(\mu
-1)k_2+s}^{\pm} & 
\end{array}
\]

These results are obtained on the basis of the isomorphism between the
algebras $\mathcal{A}_q^{-}$ (mn)$\sim {\otimes}^n\mathcal{A}_q^{-}$(m)$
\sim {\ \otimes}^m\mathcal{A}_q^{-}$(n) and the homomorphism of the
comultiplication. \bigskip

{\small \textit{Acknowledgments. } This work is supported by contract}$\Phi
-415 $ {\small with the National Fund ``Scientific Research'' of the
Bulgarian Ministry of Education and Science. }

\end{document}